\documentclass[11pt]{article}
\usepackage{amsmath,amsfonts, amssymb,graphicx,geometry,authblk,color,soul,comment,hyperref} 
\newcommand{\bff}{{\mathbf f}}
\newcommand{\bfn}{{\mathbf n}}
\newcommand{\bfzero}{{\mathbf 0}}
\newcommand{\bfK}{{\mathbf K}}
\newcommand{\bfk}{{\mathbf k}}
\newcommand{\bfu}{{\mathbf u}}
\newcommand{\bfx}{{\mathbf x}}
\newcommand{\bfe}{{\mathbf e}}
\newcommand{\bfD}{{\mathbf D}}
\newcommand{\bfy}{{\mathbf y}}

\newcommand{\bfa}{{\mathbf a}}
\newcommand{\bflambda}{{\mathbf \lambda}}

\newcommand{\bfg}{{\mathbf g}}

\usepackage[dvipsnames]{xcolor}
\usepackage{float}

\usepackage{algorithm}
\usepackage{algorithmicx}
\usepackage{algpseudocode}

\usepackage{amsthm}

\title{Accelerated Computational Micromechanics for \\ Reactive Flow in Porous Media}
\author{Mina Karimi and Kaushik Bhattacharya}
\affil{California Institute of Technology, Pasadena CA 91125}

\begin{document}
\maketitle

\begin{abstract}
    Reactive transport in permeable porous media is relevant for a variety of applications, but poses a significant challenge due to the range of length and time scales.  Multiscale methods that aim to link microstructure with the macroscopic response of geo-materials have been developed, but require the repeated solution of the small-scale problem and provide the motivation for this work. We present an efficient computational method to study fluid flow and solute transport problems in periodic porous media. Fluid flow is governed by the Stokes equation, and the solute transport is governed by the advection-diffusion equation. We follow the accelerated computational micromechanics approach that leads to an iterative computational method where each step is either local or the solution of a Poisson's equation.  This enables us to implement these methods on accelerators like graphics processing units (GPUs) and exploit their massively parallel architecture.  We verify the approach by comparing the results against established computational methods and then demonstrate the accuracy, efficacy, and performance by studying various examples. This method efficiently calculates the effective transport properties for complex pore geometries. 
\end{abstract}

\section{Introduction}

Reactive transport within permeable porous media is relevant to subsurface nuclear waste disposal \cite{spycher2003fluid}, carbon storage \cite{fan2012fully}, geothermal energy reservoirs \cite{salimzadeh2019coupled}, underground aquifers \cite{islam2016reactive},  industrial settings \cite{baqer2022review} and many other applications. However, modeling fluid and solute transport in such media presents significant challenges due to the range of interacting chemical and physical phenomena and the diversity of length and time scales.  The fluid flows through micron-size pores in a reservoir that extends over kilometers.  Chemical reactions change the solid grains, induce mineral dissolution and precipitation, evolve pore structures, and subsequently change global transport properties such as permeability and diffusivity \cite{baqer2022review}. There are empirical constitutive models that build on experimental studies \cite{baqer2022review,appelo2004geochemistry,seigneur2019reactive,witherspoon1980validity}. However, the experimental data is sparse and these models have limited fidelity and range of applicability.  Multiscale methods that aim to link microstructure with the macroscopic response of geo-materials have been developed \cite{allaire2010two,allaire2007homogenization,bear2012introduction,carbonell1983dispersion,karimi2023learning}. These methods employ homogenization techniques to derive the macroscopic equations and estimate transport properties by solving the fine-scale problem on periodic unit cells and numerical averaging \cite{allaire2010two,karimi2023learning}.   This has motivated a number of approaches to solve the fine-scale problem of fluid flow and solute transport.

One approach is to use lattice Boltzmann methods \cite{kang2006lattice,gao2017reactive} on a lattice that conforms to the microgeometry of the porous medium.  This approach is computationally expensive, and the complex and changing geometry, as well as the boundary conditions, add to the difficulty.  Another approach is to use finite element discretization of the governing continuum equations, but generating and upgrading a grid on a continually evolving complex geometry makes this extremely challenging.

This has led to a variety of methods related to the immersed boundary or fictitious domain method of Peskin \cite{peskin1972flow}.   The idea here is to extend the domain of computation to be a regular region that includes both the pores and solids, and then reintroduce the geometry through other means.  Therefore there is no need to adapt the discretization to the geometry of the porous medium.   Approaches include  Lagrange multipliers \cite{glowinski1998solution}, penalization methods \cite{kadoch2012volume}, discontinuous Galerkin \cite{chen2004pointwise,cockburn2000development}, and level set methods \cite{sethian2003level}.  It is tricky to impose general boundary conditions at the solid/fluid interface, though no-slip boundary conditions for velocity \cite{angot1999penalization}, and zero flux Neumann and Robin conditions for solute concentration \cite{kadoch2012volume} have been addressed.  We refer the reader to \cite{baqer2022review}.  In this current work, we adapt this approach of computing on an extended domain.

The fine-scale calculations are performed on a periodic domain under periodic boundary conditions, and therefore it is natural to use fast Fourier transforms (FFT) to solve the governing equations.  However, the heterogeneous coefficients pose a challenge.   Moulinec and Suquet proposed an iterative FFT method based on a Lippmann-Schwinger approximation \cite{moulinec1994fast} for the mechanical behavior of composite materials.  There are variations of this method based on an augmented Lagrangian formulation \cite{michel2000computational} and a polarization formulation \cite{monchiet_2012,eyre_1999}.   These may all be viewed as an expansion of the Neumann operator expansion \cite{moulinec_2014}.  These methods have been applied to thermoelasticity \cite{anglin2014validation}, dislocations \cite{berbenni2020fast,bertin2018fft}, shape-memory polycrystals \cite{bhattacharya2005model}, liquid crystal elastomers \cite{zhou2021accelerated}, and fluid dynamics \cite{kadoch2012volume}.  More recently, a Galerkin-Fourier approach where a Fourier basis is used for a Galerkin approximation of the underlying equations has been proposed \cite{vondvrejc2014fft} and used for calculating the effective permeability from the microstructure calculations \cite{tu2022fft}.  In this current work, we use an augmented Lagrangian approach for the fluid flow that has a variational structure and a Lippmann-Schwinger expansion for the solute transport that does not have a variational structure.

Traditionally, computational micromechanics problems have been predominantly processed on central processing units (CPUs). However, integrating CPUs and graphics processing units (GPUs) for micromechanics problems has recently gained traction to use the massively parallel architecture provided by interconnected processing units in GPU devices. GPUs can significantly accelerate Fourier transform computations due to their parallel computation \cite{bertin2018fft}.  Alternately, Hao and Bhattacharya \cite{zhou2021accelerated} proposed a formulation of the governing equations that directly exploit the structure of the GPUs.  We follow the latter approach in this work.

The goal of this work is to develop an efficient computational method to solve the problem of fluid flow and solute transport through a periodic porous medium.  The governing equations are presented in Section \ref{sec:ge}: we have the Stokes equation for the fluid flow and the advection-diffusion equation for the solute transport.  The form of these equations are chosen so that they are amenable to the calculation of macroscopic quantities in a two-scale expansion \cite{karimi2023learning}.  Our computational approach presented in Section \ref{sec:num} combines two ideas.  The first is to use an extended domain and a computational grid that does not conform to the geometry of the porous medium.  The second is to formulate the governing equations in a manner that can exploit the massively parallel architecture of GPUs.

The fluid flow is governed by the Stokes equation that may be formulated as a variational principle, and we follow the accelerated computational micromechanics approach \cite{zhou2021accelerated}, Section \ref{sec:por_form}.   We introduce the zero flow in the solid, compatibility condition, and incompressibility as constraints enforced by an augmented Lagrangian functional and solve it by the alternating direction method of multipliers (ADMM).  The result is a series of local (pointwise) calculations and a Poisson's equation, which we solve using an FFT.  This is amenable to GPU implementation, Section \ref{sec:gpu}.

The solute transport is governed by a convection-diffusion equation that is not consistent with a variational principle.   We extend the domain by introducing a small fictitious diffusivity in the solid in such a manner that it also approximates the boundary condition of the fluid-solid interface.  We can then use the 
Lippmann-Schwinger approximation \cite{moulinec1994fast}.  Again, this leads to a combination of local (pointwise) calculations and Poisson's equation, which we solve using an FFT.  This is amenable to GPU implementation, Section \ref{sec:gpu}.

We demonstrate the method in Section \ref{sec:demo} in two dimensions using two geometries.  The first is a model geometry that has been widely used as a benchmark, and the second is obtained by sampling a real material.  The proposed computational approach in both the fluid flow and solute transport is iterative, and the convergence depends on the choice of a few (numerical) parameters.  We study this and propose an efficient approach to choosing these parameters.  We then verify the numerical approach by comparing the solutions obtained by the proposed method with that of a finite element method on a grid that conforms to the microstructure.  Finally, we demonstrate that the approach has good parallel performance, and importantly, the GPU implementation is significantly faster than a CPU implementation using a variety of methods.

\section{Reactive flow in periodic porous medium} \label{sec:ge}

\begin{figure}
\centering
    \includegraphics[width=4in]{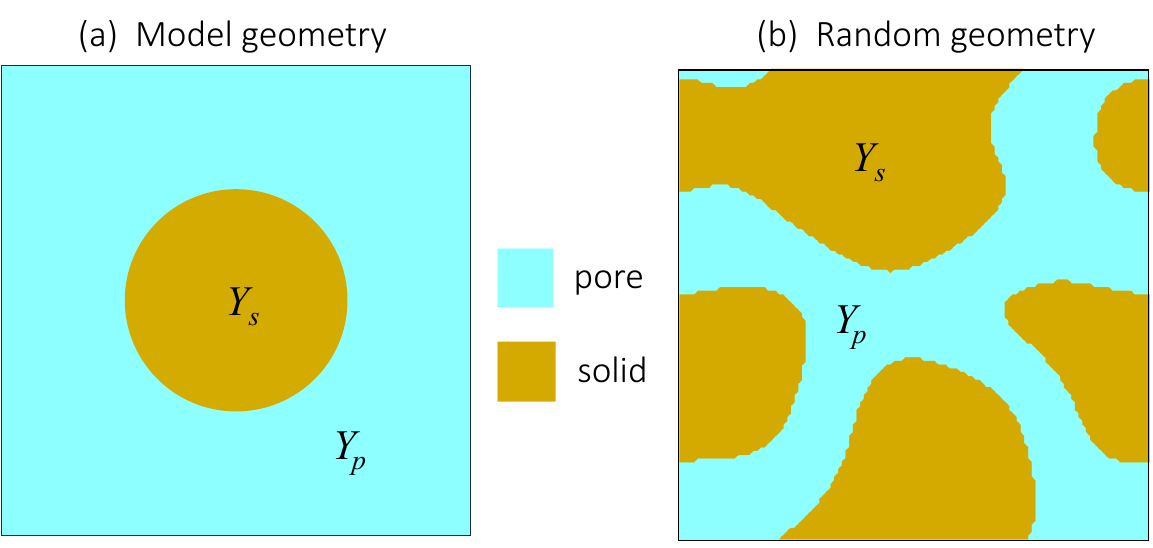}
    \caption{Two examples of the unit cell.  (a) Model geometry that is widely used as a benchmark, and (b) Random geometry.}
    \label{fig:geom}
\end{figure}

We consider a periodic porous medium with a unit cell $Y$ consisting of a porous region $Y_p$ and solid region $Y_s$ with $Y = Y_p\cup Y_s$.  We choose $Y$ to be the unit cube by suitable choice of length scales.  See Figure \ref{fig:geom} for two examples.  We consider a flow through the periodic medium driven by an average or overall pressure gradient $-\bfg_p$.  We assume that the fluid carries a solute and the transport is driven by an average or overall chemical composition gradient $\bfg_\chi$.   The fluid flow is governed by the Stokes equation with the no-slip boundary condition at the solid/fluid interface, and the solute transport is governed by the advection-diffusion equation with a no flux boundary condition at the interface:
\begin{align} \label{eq: unit cell}
    \begin{cases}
            -\nabla p + \nu \Delta \bfu + \bfg_p = \bfzero, \ \  \nabla \cdot \bfu = 0 \quad & \text{in} \ Y_p \\
	    \bfu = \bfzero & \text{in} \ Y_s \\
	    \text{Pe} \ \bfu\cdot \left( \nabla \chi \right)  - \Delta \chi  = \
	    \text{Pe} (\bar \bfu - \ \bfu)\cdot \bfg_\chi \quad & \text{in} \ Y_p\\
       - \left( \nabla \chi \right)\cdot {\bfn} =  \bfg_\chi \cdot {\bfn} & \text{on} \ \partial Y_p	    
    \end{cases}
\end{align}
where $\bfu$ denotes the fluid velocity and $p$ the excess pressure (so that pressure is $-\bfg_p \cdot \bfx + p$) and $\chi$ denotes the excess chemical composition (so that the composition is $ \bfg_\chi \cdot \bfx + \chi$). $\nu$ is the viscosity and $\text{Pe}$ the Pecl\'et number.  $\bar \bfu = 1/|Y_p \int_{Y_p} \bfu dy$ is the pore-averaged velocity.

The linearity of Stokes flow (\ref{eq: unit cell})$_{1,2}$ means that we can find the solution as $\bfu = \sum_{i=1}^d (g_p)_i \bfu^i$ where $\bfu^i$ solves the problem with the average pressure gradient set to $\hat \bfe_i, i = 1, \dots, d$, an orthonormal basis aligned with the unit cell.  Similarly, the linearity of the transport equation (\ref{eq: unit cell})$_{3,4}$ implies that $\chi = \sum_{i=1}^d (g_\chi)_i \chi^i$ where $\chi^i$ solves the problem with average chemical gradient set to $\hat \bfe^i$.

This system of equations can be used for the microscopic scale in two-scale expansion to obtain the macroscopic permeability $\bfK^*$ and diffusivity $\bfD^*$:
\begin{align}
\bfK^*_{ij} &= \int_{Y_p} \nabla \bfu^i \cdot \nabla \bfu^j \ dy, \label{eq:K}\\
\bfD^*_{ij} &= \int_{Y_p} \bfe^i\cdot \bfe^j \ dy + \text{Pe} \int_{Y_p} (\overline \bfu^i -  \bfu^i) \chi^j \ dy 
    +  \int_{Y_p} \nabla \chi^j\cdot \bfe^i \ dy . \label{eq:D}
 \end{align}
where $\bfu^i, p^i$ and $\chi^i$ solve (\ref{eq: unit cell}) with $\bfg_p = \hat \bfe^i$ and  $\bfg_\chi = \hat \bfe^i$ respectively.  Note that the macroscopic diffusivity may not be symmetric due to the influence of the flow at the microscopic scale.

When one has deposition or dissolution, the solid region $Y_s$ evolves slowly over time, and one has to repeatedly solve the unit cell problem (\ref{eq: unit cell}).  This motivates the development of an efficient numerical method.  We refer the reader to Karimi and Bhattacharya \cite{karimi_2023} for a discussion of the two-scale expansion in the presence of dissolution/deposition.

\section{Numerical method and implementation} \label{sec:num}

We seek an efficient and fast numerical method to solve the governing equations (\ref{eq: unit cell}).   We do so with two ideas motivated by the accelerated computational micromechanics (ACM) approach \cite{zhou2021accelerated}.   The first is to formulate the governing equations (\ref{eq: unit cell}) on the entire unit, and use a fast Fourier transform based numerical method.   The second is to implement the method efficiently on accelerators or graphical processing units (GPUs). 

We introduce the indicator or characteristic function $H$ for the solid region $Y_s$:
\begin{equation}\label{eq: Heaviside}
    H(\bfy)=
    \begin{cases}
        1 & \text{if } \ \bfy \in Y_s,\\
        0 & \text{if } \ \bfy \in Y_p.
    \end{cases}
\end{equation}

\subsection{Porous media flow} \label{sec:por_form}

We begin with the flow equations (\ref{eq: unit cell})$_{1,2}$.  We can rewrite it as a constrained variational principle,
\begin{equation}
\min_{\nabla \cdot \bfu = 0,\ \bfu|_{Y_s} = 0} \ \int_Y \left( \frac{\nu}{2} |\nabla \bfu|^2 - \bfu \cdot \bfg_p \right) \ dy.
\end{equation}
Therefore it is readily amenable to ACM.    We treat the constraints using an augmented Lagrangian and solve the problem using the alternating direction method of multipliers (ADMM) \cite{glowinski1989augmented}.  

Consider the augmented Lagrangian functional 
\begin{multline}  \label{eq:Lagrangian with auxiliary var}
   \int_Y \frac{\nu}{2} |\nabla \bfu|^2 - \bfu\cdot \bfg_p
    - \left( q(\nabla\cdot\bfu) - \frac{\beta}{2}|\nabla\cdot\bfu|^2 \right) 
    + \left(\bfa\cdot\left(\bfu - \Tilde{\bfu} \right) + \frac{b}{2}|\bfu- \Tilde{\bfu}|^2\right) + 
    H\left( \bflambda\cdot \Tilde{\bfu} + \frac{\alpha}{2}| \Tilde{\bfu}|^2 \right) 
\end{multline}
where we have introduced an auxiliary variable $\tilde \bfu$, and Lagrange multipliers $q, \bfa, \bflambda$ and penalty coefficients $\beta>0, b>0, \alpha>0$.  The two terms in the first parenthesis are introduced to ensure incompressibility $\nabla \cdot \bfu = 0$, those in the second term enforce $\bfu = \tilde \bfu$ and those in the final enforce $\tilde \bfu = \bf0$ in $Y_s$.

We solve the constrained minimization problem using the iterative alternating direction method of multipliers (ADMM) ~\cite{glowinski1989augmented, boyd2011distributed}.  Set $X^n=\{\bfu^n, \tilde \bfu^n, q^n, \bfa^n, \bflambda^n \}$.  We initialize the iteration with the guess $X^0$.  At the $n^\text{th}$ step, we obtain $X^{n+1}$ from $X^n$ as follows:
\begin{itemize}
\item[]\textit{Step 1. Stationarity in $\bfu$.} Find $\bfu^{n+1}$ by solving the partial differential equation
\begin{align} \label{eq:step1por}
    -\nu \nabla\cdot(\nabla\bfu^{n+1})-\hat\bfe + \nabla q^n - \beta\nabla(\nabla\cdot\bfu^{n+1}) + \bfa^n + b\left( \bfu^{n+1} - \Tilde{\bfu}^n \right)  = 0;
\end{align}
\item[]\textit{Step 2. Stationarity in $\tilde \bfu$.} Update $\Tilde{\bfu}$ by solving at each point in the unit cell
\begin{equation} \label{eqn:u_tilde}
    -\bfa^n - b\left( \bfu^{n+1} - \Tilde{\bfu}^{n+1} \right) + H\bflambda^n  +\alpha H\Tilde{\bfu}^{n+1} = 0;
\end{equation}
\item[]\textit{Step 3. Lagrange multiplier update.} Set at each point in the unit cell
\begin{equation}\label{eqn:Lagrane-multipliers}
\begin{aligned} 
    q^{n+1} &= q^n - \beta\left( \nabla\cdot\bfu^{n+1} \right) \ , & 
    \bfa^{n+1} &= \bfa^n + b \left( \bfu^{n+1} - \Tilde{\bfu}^{n+1} \right) \\
    \bflambda^{n+1} &= \bflambda^n + \alpha \left( H\Tilde{\bfu}^{n+1} \right). 
\end{aligned}
\end{equation}
\end{itemize}
We iterate till we converge according to the following criteria:
\begin{itemize}
\item[]\textit{Step 4. Check for convergence} 
\begin{equation}\label{eq:constrains}
\begin{aligned}
    r_{p1} &= ||H\Tilde{\bfu}^{n+1}||_{L^2}\leq r_{p1}^\text{tol}, & r_{d1} &= \alpha|| H\left(\Tilde{\bfu}^{n+1} - \Tilde{\bfu}^n\right) ||_{L^2}\leq r_{d1}^\text{tol}, \\ 
    r_{p2} &=||\nabla\cdot\bfu^{n+1}||_{L^2}\leq r_{p2}^\text{tol}, & r_{d2} &= \beta|| \nabla\cdot\bfu^{n+1} - \nabla\cdot\bfu^n ||_{L^2}\leq r_{d2}^\text{tol}, \\
    r_{p3} &= || \bfu^{n+1} - \Tilde{\bfu}^{n+1} ||_{L^2}\leq r_{p3}^\text{tol}, &r_{d3} &= b|| \bfu^{n+1} - \bfu^n ||_{L^2}\leq r_{d3}^\text{tol} \
\end{aligned}
\end{equation}
\end{itemize}
for given primal and dual tolerance values $r_{p1}^\text{tol}, ~ r_{d1}^\text{tol}, ~ r_{p2}^\text{tol}, ~ r_{d2}^\text{tol}, r_{p3}^\text{tol}$ and $r_{d3}^\text{tol}$  respectively. The tolerance values are defined using absolute and relative criteria, with specific values of relative and absolute tolerance depending on the application \cite{boyd2011distributed}. The tolerance values are normalized by the number of degrees of freedom; as an example, the tolerance values for the first constraint $r_{p1}^{\text{tol}}$, and $r_{d1}^{\text{tol}}$ are defined as:
\begin{align}
    r_{p1}^{\text{tol}} &= \sqrt{n}\epsilon_{\text{abs}} + \epsilon_{\text{rel}}\max \{ ||H\tilde{\bfu}^{n+1}||_{L^2}, ||\lambda^{n+1}||_{L^2} \} \\
    r_{d1}^{\text{tol}} &= \sqrt{n}\epsilon_{\text{abs}} + \epsilon_{\text{rel}} ||\lambda^{n+1}||_{L^2}
\end{align}
where $n$ is the number of degrees of freedom for $\Tilde{\bfu}$, and $\epsilon_{\text{abs}}$, and $\epsilon_{\text{rel}}$ are the positive absolute and relative tolerances, respectively. The relative term can regulate the tolerance value and convergence speed. Specifically, in the case of time-dependent geometries, the relative term prevents fast or slow convergence. In the following, we consider the tolerance $\epsilon = \epsilon_{\text{abs}} = \epsilon_{\text{rel}}$.

The only differential equation is Step 1, but this is linear, involving only constant coefficients.  So, it can be solved using the fast Fourier transform (FFT) and its inverse as follows
\begin{equation}\label{eq:FFT step1}
    \nu |\bfk|^2 \hat{\bfu}^{n+1} - \hat{\bfg_p} + \bfk \hat{q}^n + \beta \bfk (\bfk \cdot \hat{\bfu}^{n+1}) + \hat{\bfa}^n + b (\hat{\bfu}^{n+1} - \hat{\Tilde{\bfu}}^n) = \bfzero
\end{equation}
Steps 2 and 3 are locally solved at each point.  For these reasons, this algorithm can be implemented efficiently in parallel using accelerators like GPUs.

The method convergences for a sufficiently large value of the penalty parameters $\alpha, \beta, b$, but the rate of convergence depends on the values  \cite{boyd2011distributed}.  We study this in Section \ref{sec:demo}.

\subsection{Solute transport}\label{sec:concentration iterative algorithm}

We now turn to solute transport, (\ref{eq: unit cell})$_{3,4}$.  These equations are not associated with a variational principle, and so we proceed differently.    We introduce a small constant (fictitious diffusivity) $\eta$ to be specified, and consider the equation 
\begin{equation} \label{eq: penalization}
   \text{Pe} \left(1-H \right) \bfu\cdot\left( \nabla\chi + \bfg_\chi \right) - \nabla\cdot\left[ \left(\left(1-H \right) + \eta H \right)\left(\nabla\chi + \bfg_\chi \right)\right] = \text{Pe} \left( 1-H\right)\bar \bfu \cdot \bfg_\chi
\end{equation}
on $Y$ for periodic $\chi$ where $H$ is the characteristic function of $Y_s$ (\ref{eq: Heaviside}).  Note that $\nabla \chi$ may suffer a discontinuity across the pore-solid interface (due to the discontinuous coefficient in the second order term); however, the weak form leads to the jump condition
\begin{equation} \label{eq:jump}
\text{Pe} (\left. \nabla \chi \right|_p + \bfg_\chi) \cdot \bfn = \eta (\left. \nabla \chi \right|_s + \bfg_\chi) \cdot \bfn
\end{equation}
at any point on the interface between the pore and solid where $\left. \nabla \chi \right|_p$ and $\left. \nabla \chi \right|_s$ are limiting values of the gradient from the pore and solid, respectively.  Thus, the solution of (\ref{eq: penalization}) provides an approximation of (\ref{eq: unit cell})$_{3,4}$ as $\eta \to 0$.  

Note that the highest order term in (\ref{eq: penalization}) has a non-constant coefficient.  Therefore, a Fourier transform of this equation would lead to convolutions.  So, we follow Moulinec and Suquet\cite{moulinec_1994} (and many others since) to rewrite it as
\begin{equation}
    - \mathcal{A}_0 \Delta \chi + \mathcal{B}_0 \cdot \left( \nabla \chi +\bfg_\chi \right) = 
    \mathcal{F} + \nabla\cdot \left[ (\mathcal{A} - \mathcal{A}_0)  (\nabla\chi + \bfg_\chi ) \right] - (\mathcal{B} - \mathcal{B}_0) \cdot \left( \nabla \chi + \bfg_\chi \right) \end{equation}
where $\mathcal{A} = (1-H) + \eta H, \ \mathcal{B} = \text{Pe}(1-H)\bfu , \ \mathcal{F} = \text{Pe} (1-H) \bar \bfu \cdot \bfg_{\chi}$ and $\mathcal{A}_0 > 0, \mathcal{B}_0$ are uniform ``comparison medium values'' to be specified.
We can solve this equation by iteration 
\begin{equation} \label{eq:iter}
- \mathcal{A}_0 \Delta \chi^{n+1} + \mathcal{B}_0 \cdot \left( \nabla \chi^{n+1} +\bfg_\chi \right) = F^n
\end{equation}
where 
\begin{equation}
F^n = \mathcal{F} + \nabla\cdot \left[ (\mathcal{A} - \mathcal{A}_0)  (\nabla\chi^n + \bfg_\chi ) \right] - (\mathcal{B} - \mathcal{B}_0) \cdot \left( \nabla \chi^n + \bfg_\chi \right) .
\end{equation}

We use the following iterative algorithm.  We initialize the calculation with an initial guess $\chi^0, \nabla \chi^0$.  At the $n^\text{th}$ step, we obtain $\chi^{n+1}, \nabla \chi^{n+1}|$ as follows:
\begin{itemize}
\item[]\textit{Step 1. Compute the residual.} Compute $F^n$ in real space;
\item[]\textit{Step 2. Update concentration.} Solve (\ref{eq:iter}) for $\chi^{n+1}$.
\end{itemize}
We iterate till we converge according to the following criteria:
\begin{itemize}
\item[]\textit{Step 3. Check for convergence} 
\begin{equation}\label{eqn:tol_solute}
\begin{aligned}
    r_1 &= ||\chi^{n+1}-\chi^n||_{L^2} \leq r_1^{\text{tol}}, & r_2 = ||\nabla\chi^{n+1}-\nabla\chi^n||_{L^2} \leq r_2^{\text{tol}}
\end{aligned}
\end{equation}
\end{itemize}
where $r_1^{\text{tol}}$, and $r_2^{\text{tol}}$ tolerances are defined as $\sqrt{n}\epsilon$, and $n$ is the dimension of $\chi$, and $\nabla\chi$ fields. 

The only differential equation is in Step 2, but this is linear with constant coefficients.  So, we can solve it efficiently in Fourier space as
\begin{equation}\label{eqn:fft_solute}
    \hat{\chi}^{n+1} = \frac{\hat{F}^n - \mathcal{B}_0\cdot\hat{\bfg_\chi}}{\mathcal{B}_0\cdot i \bfk +\mathcal{A}_0 |\bfk|^2}.  
\end{equation}
There is, however, one aspect that requires some care.  Note that we expect the concentration gradient $\nabla \chi$ to be discontinuous across the solid-fluid interface (cf. (\ref{eq:jump})).  This can lead to spurious oscillations when we use discrete Fourier transforms in computations \cite{berbenni2014numerical}.   To overcome this problem, we replace the discrete Fourier transform of the derivatives with the discrete Fourier transform of central difference approximation:
\begin{align}  \label{eq:fft_app}
    \mathcal{F}(\chi_{,j}) &=ik_j\mathcal{F}(\chi) \approx \frac{ \mathcal{F}\left[\chi(x+he_j)\right] - \mathcal{F}\left[\chi(x-he_j)\right]}{2h} = i\mathcal{F}(\chi) \frac{\sin(hk_j)}{h} \\ \nonumber
    \mathcal{F}(\chi_{,jj}) &= -|k_j|^2\mathcal{F}(\chi)  \approx \frac{\mathcal{F}\left[\chi(x+he_j)\right] - \mathcal{F}\left[\chi(x-he_j)\right] - 2\mathcal{F}(\chi)}{h^2}
    -\mathcal{F}(\chi) \frac{4\sin^2(\frac{hk_j}{2})}{h^2}.
\end{align}
This acts as a high-frequency filter that is known to suppress spurious oscillations \cite{berbenni2014numerical, lebensohn2016numerical, zhou2021accelerated}.

Step 1 is a local update and therefore the iterative algorithm can be efficiently implemented on GPUs.

The method converges for a sufficiently large value of the comparison medium values  $\mathcal{A}_0, \mathcal{B}_0$, but the rate of convergence depends on the values \cite{moulinec_2018}.  We study this in Section \ref{sec:demo}.

\subsection{GPU implementation} \label{sec:gpu}

\begin{algorithm}[t]
\caption{GPU Implementation }\label{alg:parallel}
\begin{algorithmic}
\State Given the microstructure geometry $H$; 
\State \textit{Initialize}.
\State Place $\epsilon$, $\alpha$, $\beta$, $b$ in the constant cache;
\State Place $\mathcal{A}_0$, $\mathcal{B}_0$ in the global memory;\\
    \State *\underline{Problem 1: Porous media flow}.
    \State Initialize $\bfu^0$, $\Tilde{\bfu}^0$, $q^0$, $\bfa^0$, $\bflambda^0$;
    \State Place $\bfu^0$, $\Tilde{\bfu}^0$, $q^0$, $\bfa^0$, $\bflambda^0$ in the global memory;
    \While {$\left( r_{p1}>r_{p1}^{\text{tol}}~\text{or}~r_{d1}>r_{d1}^{\text{tol}} \right)$ and $\left( r_{p2}>r_{p2}^{\text{tol}}~\text{or}~r_{d2}>r_{d2}^{\text{tol}} \right)$ and $\left( r_{p3}>r_{p3}^{\text{tol}}~\text{or}~r_{d3}>r_{d3}^{\text{tol}} \right)$}
    \State \textit{Step 1: Stationary in} $\bfu$. 
    \State \ \ \ \ \ - FFT $\hat{\bfu}^{n+1}$, $\hat{\Tilde{\bfu}}^{n}$, $\hat{q}^{n}$, $\hat{\bfa}^{n}$; \Comment{ /* cuFFT */}
    \State \ \ \ \ \ - Compute $\hat{{\bfu}}^{n+1}$ from (\ref{eq:FFT step1}); 
    \State \ \ \ \ \ - iFFT $\hat{{\bfu}}^{n+1}$; \Comment{/* cuFFT */}
    \State \textit{Step 2: Stationary in} $\Tilde{\bfu}$. Compute $\Tilde{\bfu}^{n+1}$ from (\ref{eqn:u_tilde});
    \State \textit{Step 3: Lagrange multiplier update}. Compute $q^{n+1}$, $\bfa^{n+1}$, $\bflambda^{n+1}$ from (\ref{eqn:Lagrane-multipliers});
    \State \textit{Step 4: Check for convergence}. Compute $r_{p1}$, $r_{d1}$, $r_{p2}$, $r_{d2}$, $r_{p3}$, $r_{d3}$; \Comment{/* cuBlas */}
    \EndWhile \\
    \State *\underline{Problem 2: Solute transport}.
    \State Initialize $\chi^0$, $\nabla\chi^0$;
    \State Place $\chi^0$, $\nabla\chi^0$ in the global memory;
    \While {$\left( r_1 > \sqrt{n}\epsilon \right)$ and $\left( r_2 >\epsilon \right)$}
    \State \textit{Step 1: Compute the residual} Compute $F^n$.
    \State \textit{Step 2: Update the concentration}
    \State \ \ \ \ \ - FFT $\hat{F}^{n+1}$, $\hat{\chi}^{n+1}$; \Comment{/* cuFFT */}
    \State\ \ \ \ \  - Compute $\hat{\chi}^{n+1}$ from (\ref{eqn:fft_solute});
    \State\ \ \ \ \  - iFFT $\chi^{n+1}$; \Comment{/* cuFFT */}
    \State \textit{Step 3: Check convergence}. Compute $r_1$, $r_2$ from (\ref{eqn:tol_solute}); \Comment{/* cuBlas */}
    \EndWhile
\end{algorithmic}

\end{algorithm}

Algorithm \ref{alg:parallel} summarizes the algorithms described in the previous sections, and outlines the implementation on GPUs.
We refer the reader to \cite{kirk_2016} for details of the GPU architecture.   In this work, we used NVIDIA P100 GPUs equipped with a compute unified device architecture (CUDA) platform with the C++ programming language. The calculation is initiated by a main processor (CPU).   Data transfer between a CPU and GPU is slow, and therefore we limit data exchange only to initialization and the recovery of results for storage.  The rest of the calculation is conducted entirely in the GPU.   All the local calculations are trivially parallel and, therefore, solved independently at each point in a separate thread.  The FFT and iFFT are solved using a standard library ({\tt cuFFT} in CUDA).  The check for convergence requires the computation of norms and this is also done using a standard linear algebra library ({\tt cuBlas} in CUDA).

\section{Convergence, verification and performance} \label{sec:demo}
We now demonstrate the performance and convergence of the iterative algorithms, considering unit cells with a model pore geometry that has been widely used as a benchmark as well as a realistic random pore geometry, see Figure \ref{fig:geom}. 

\subsection{Porous media flow}

\paragraph{Penalty parameters and convergence.}
\begin{figure}
\centering
    \includegraphics[width=6in]{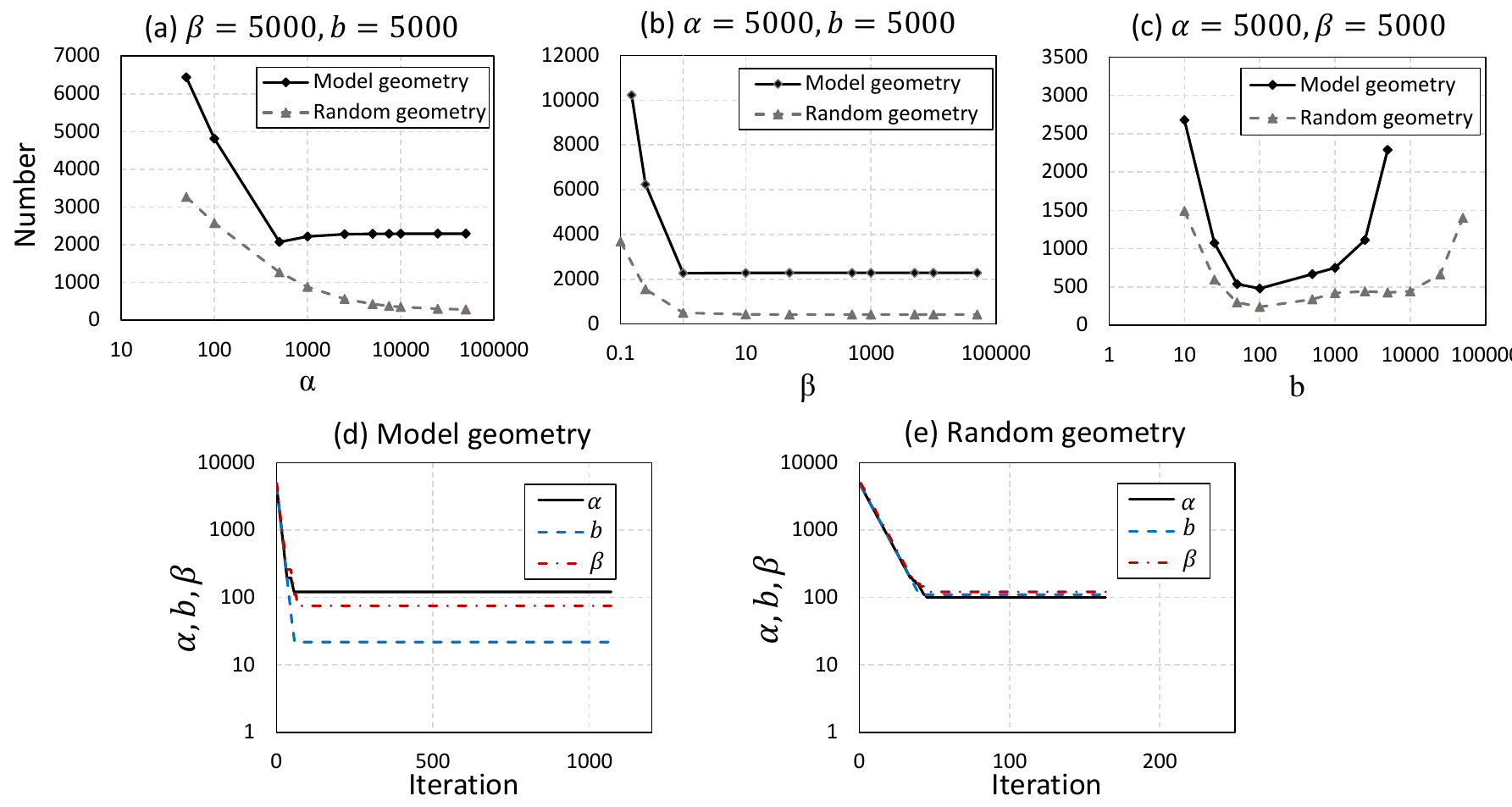}
    \caption{Effect of penalty parameters on convergence of the Stokes flow algorithm.  (a-c) The number of iterations required for convergence for various choices of fixed penalty parameters for simulations in the model geometry.  (d,e) The evolution of penalty parameters with iteration according to (\ref{eq: varying algorithm}) for the model (d) and random (e) geometries.   Tolerance $\epsilon=10^{-5}$ in these calculations.}
    \label{fig:pen_par}
\end{figure}

It is known that the ADMM method described in Section \ref{sec:por_form} and Algorithm \ref{alg:parallel} converges for sufficiently large values of the penalty parameters $\alpha$, $\beta$, and $b$ ~\cite{boyd2011distributed}. However, the rate of convergence is highly sensitive to the values of penalty parameters.   This is shown in Figure \ref{fig:pen_par}(a-c).  Typically, larger (respectively smaller) values of the penalty parameter result in a faster drop in the value of primal error $r_p$ (respectively dual error $r_d$) and a slower decrease in dual error $r_d$ (respectively primal error $r_p$).   This can result in slow convergence.  

To address this limitation, an extension of this approach has been proposed ~\cite{boyd2011distributed, zhou2021accelerated}, where penalty parameters are varied as iteration proceeds.  This approach adjusts the value of $\psi = \{\alpha, \beta, b\}$ within the iterations to improve the convergence rate as follows:\begin{equation}\label{eq: varying algorithm}
    \psi_k^{n+1} = 
    \begin{cases}
        \phi_k\psi_k^n, & \text{if} ~~~ r_p/r_d>\tau_k \\
        \max\{\psi_k^n/\phi_k, ~ {\psi_k}_{\min}\}, & \text{if} ~~~ r_d/r_p>\tau_k \\
        \psi_k^n, & \text{if} ~~~ else
    \end{cases}
\end{equation}
Figure \ref{fig:pen_par}(d,e) shows the variation of penalty parameters with iteration using the varying penalty parameter approach  (\ref{eq: varying algorithm}). In this simulation, we set the $\phi_k$, and $\tau_k$ values for penalty parameters to be $\phi=\{ 1.1, 1.1, 1.1 \}$, and $\tau = \{ 20, 10, 30 \}$\footnote{We use $\tau = \{ 5, 10, 5 \}$ for the random geometry below.}.   Note that the number of iterations required for convergence is significantly smaller than those for fixed parameters.  We show in Supplementary Materials, Figure \ref{fig:prim_dual}, that the residuals do not decrease monotonically as we vary the penalty parameters as the iteration proceeds with varying penalty parameters.

\paragraph{Verification.}  We now turn to the verification of the numerical method.  We solve the flow problem with three different approaches

\begin{itemize}
\item {\it Fluid domain FEM.}  We solve the original Stokes flow equation (\ref{eq: unit cell})$_{1}$ on the fluid domain $Y_p$ subject to periodic boundary conditions on the fluid boundary of the unit cell $\partial Y_p \cap \partial Y$ and to Dirichlet boundary condition $\bfu = \bf0$ on the solid-fluid interface $Y_p \cap Y_s$.  We use an unstructured mesh to resolve the equations on $Y_p$.  
\item {\it Extended domain FEM.}  We solve the ADMM-based formulation on the extended domain as described in Section \ref{sec:por_form} by using a structured finite element discretization of (\ref{eq:step1por}) on the unit cell $Y$.
\item {\it Proposed method.} We solve the ADMM-based formulation on the extended domain using  FFT and GPUs as described in Section \ref{sec:gpu}.
\end{itemize}

\begin{figure}
\centering
    \includegraphics[width=6in]{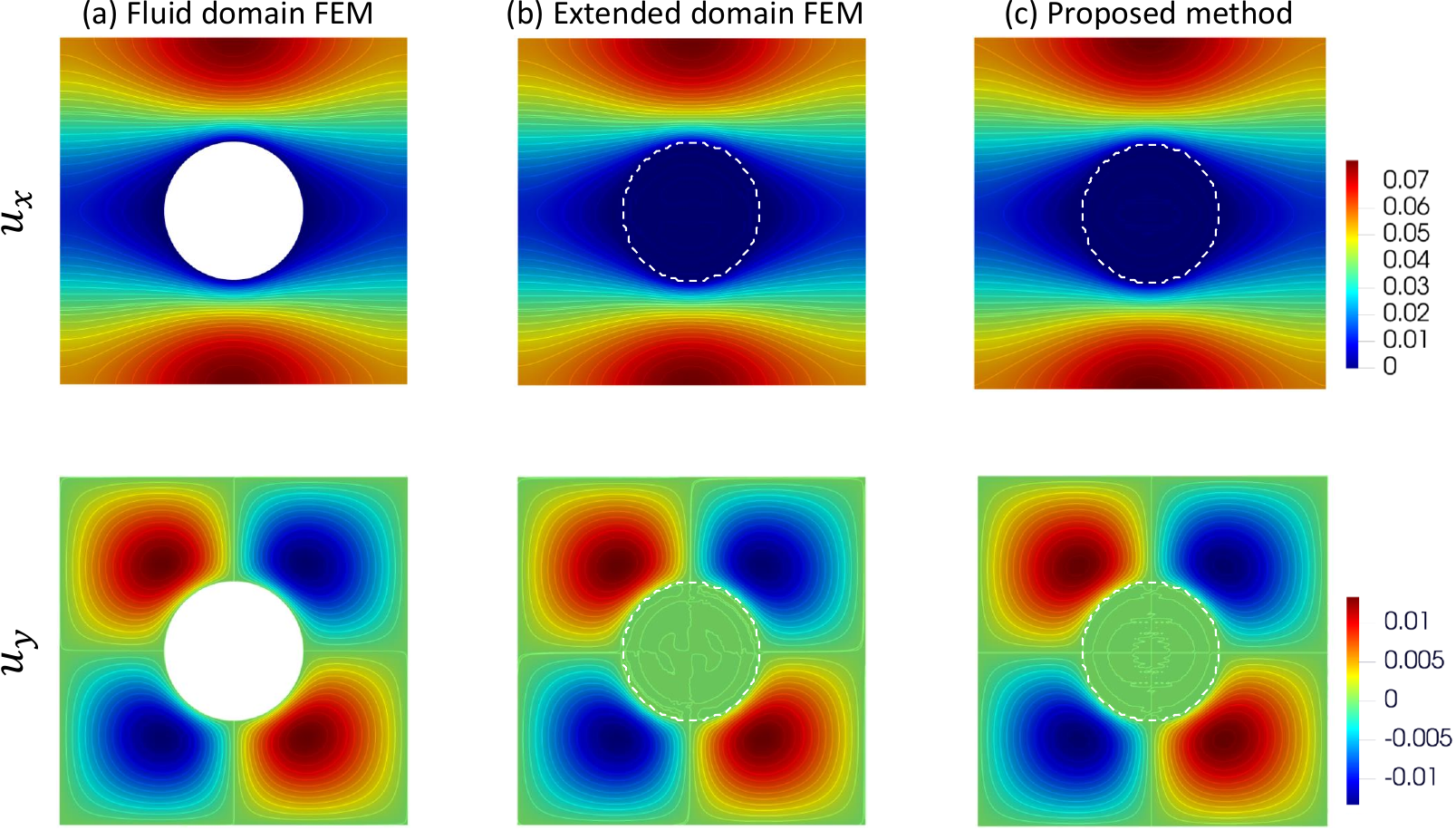}
    \caption{Verification of flow in the model geometry by comparing the solution using the three methods with $\bfg_p = [1 ~ 0]$ and $\epsilon = 10^{-5}$. The white dashed line indicates the solid/fluid interface.}
    \label{fig:model_vel}
\end{figure}

\begin{figure}
\centering
    \includegraphics[width=4.5in]{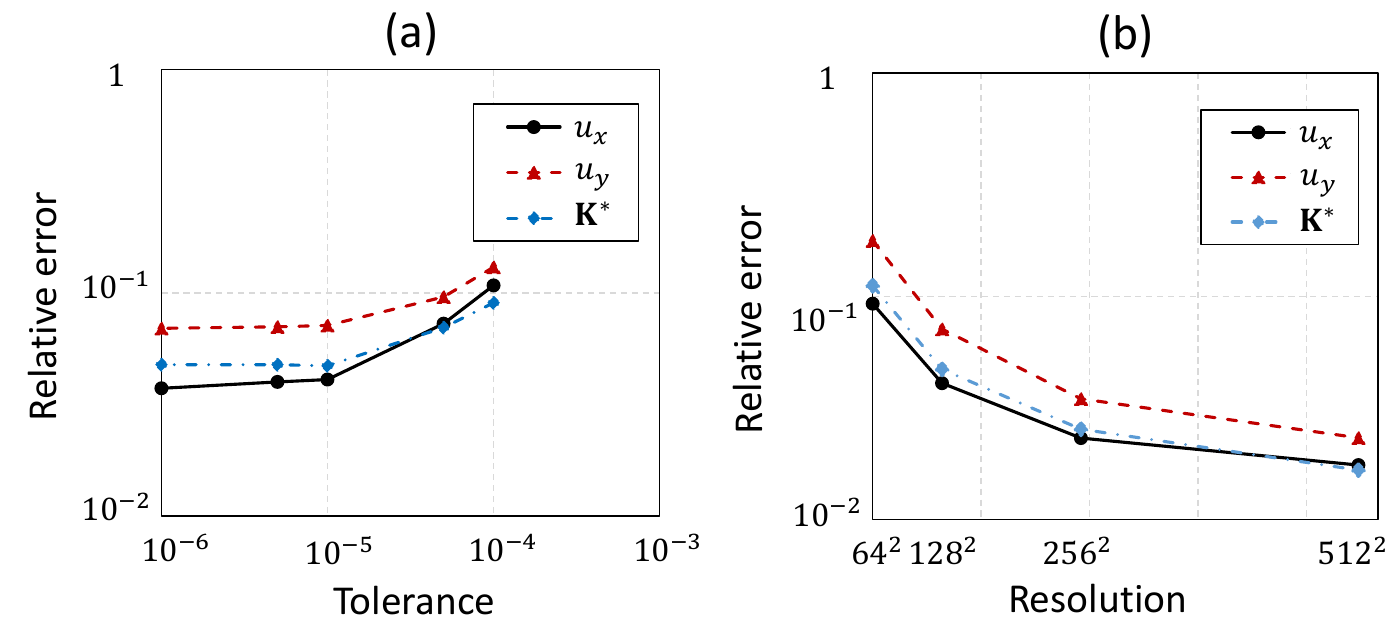}
    \caption{The relative errors in velocity fields and the permeability computed using the proposed method relative to the fluid domain FEM in the model geometry for various  (a) tolerances and (b) resolution.    $\bfg_p = [1 ~ 0], \epsilon=10^{-5}$.} 
    \label{fig:model_err}
\end{figure}

Figure \ref{fig:model_vel} presents a comparison of velocity profiles obtained using three methods for the model geometry.  The fluid domain FEM used an unstructured mesh consisting of 100,856 triangular elements, while the other two methods used a structured discretization with a resolution of $128\times 128$.  We observe that the constrained formulation in both the extended domain FEM and proposed FFT-based method results in zero velocity in the solid region as required and that all results of all three methods match very well in the fluid or pore region.  The relative  $L_2$ norm (least squares) over the entire unit cell\footnote{We have verified that results are similar if we use the relative  $L_2$ norm (least squares) over the pore space} if we use the  of the velocity components computed using the proposed method relative to the fluid domain FEM is shown in Figure \ref{fig:model_err} for various values of tolerance ($\epsilon$) and resolution. The relative error in the computed permeability (\ref{eq:K}) is also shown.  We find that the errors are small, and also decrease with lower tolerance and finer resolution.

\begin{figure}
\centering
    \includegraphics[width=3.8in]{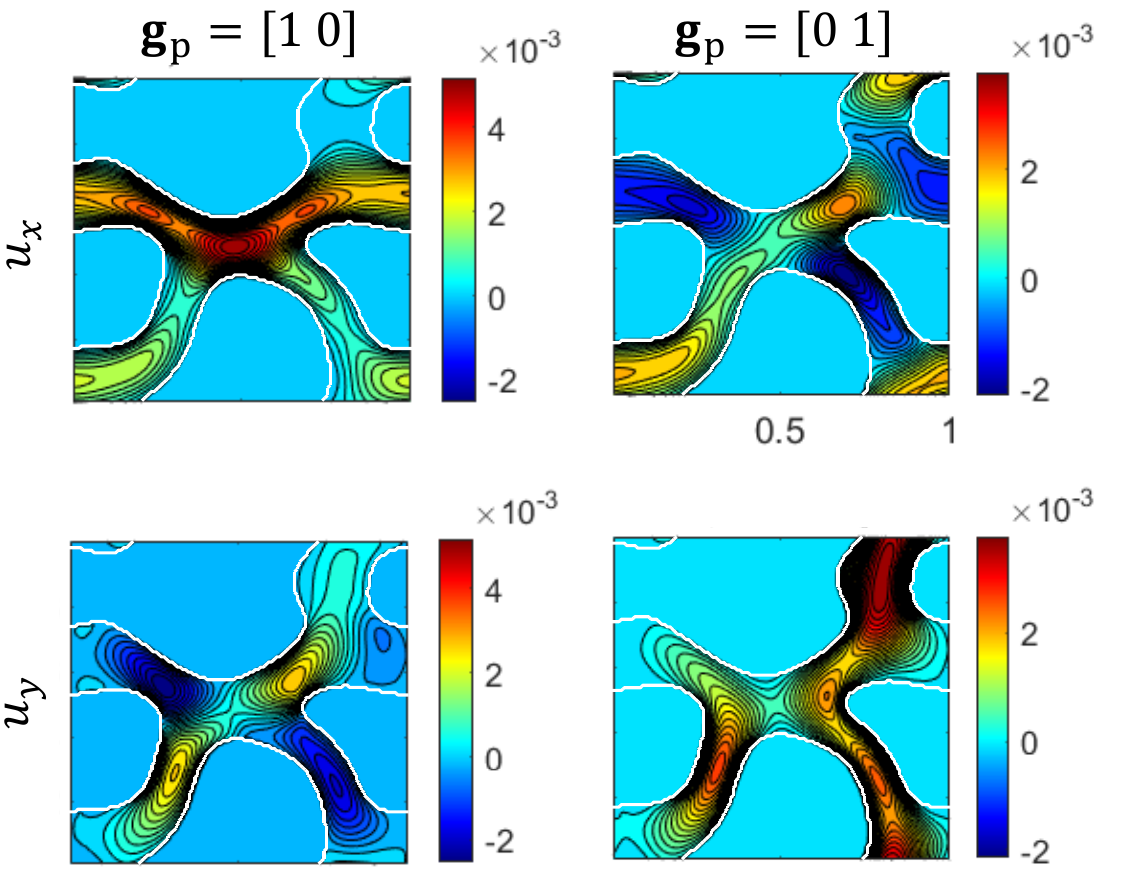}
    \caption{Anisotropy of flow induced by the geometry.  $\epsilon=10^{-5}$ in these calculations. }
    \label{fig:vel_random}
\end{figure}
We find a similar agreement between the extended domain FEM and the proposed FFT-based method in the random geometry (see Supplementary Materials, Figures \ref{fig:random_vel1}, \ref{fig:random_vel2} and \ref{fig:random_err}).  We have not performed the fluid domain FEM in this case.

We conclude this by noting in Figure \ref{fig:vel_random} the anisotropy in the flow induced by the anisotropy of the microstructure.   We see that the random microstructure leads to different flow fields when $\bfg_p = [1 \ 0]$ and $\bfg_p = [0 \ 1]$.  This leads to anisotropy in the permeability $\bfK^*$ (and anisotropy and asymmetry in the diffusion coefficient $\bfD^*$, as we see later).


\subsection{Solute Transport}

\paragraph{Comparison medium and convergence.}
\begin{figure}
\centering
    \includegraphics[width=4.5in]{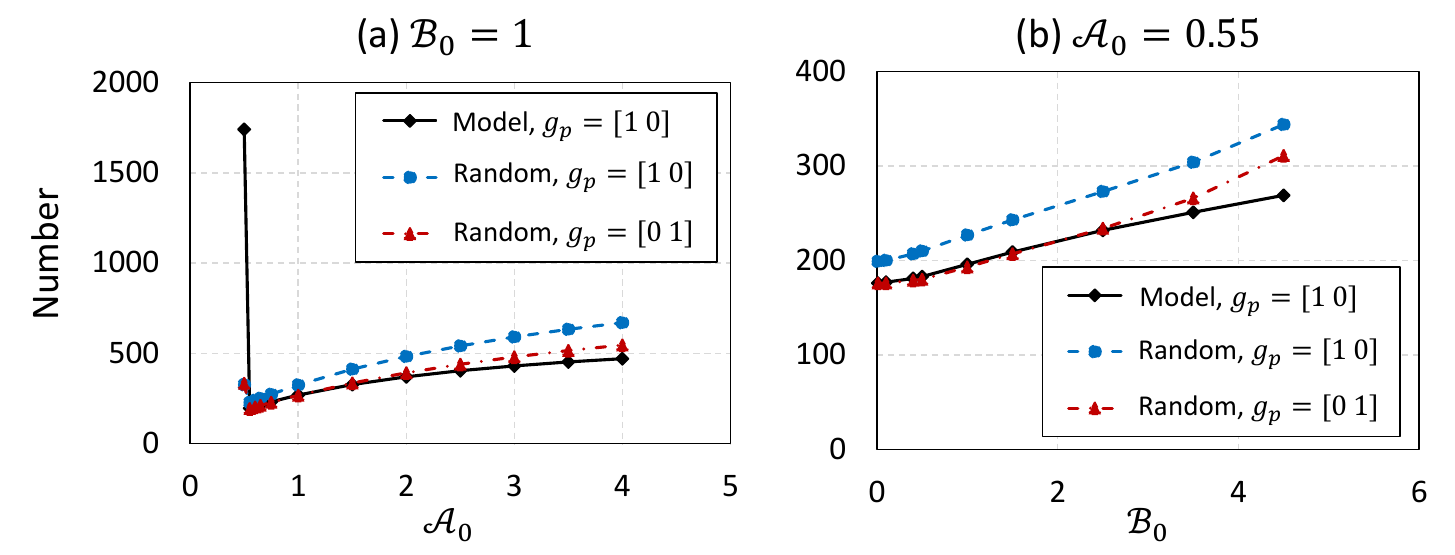}
    \caption{The effect of comparison medium values on the convergence of solute transport in the model geometry.  (a) The number of iterations required convergence vs.\  ${\mathcal A}_0$ with fixed ${\mathcal B}_0$. (b) The number of iterations required convergence vs.\  ${\mathcal B}_0$ with fixed ${\mathcal A}_0$.  $\bfg_\chi = [1 ~ 0]$, $ \eta = 0.01$, $\epsilon=10^{-5}$.} 
    \label{fig:comp_med}
\end{figure}

It is known that the method proposed in Section \ref{sec:concentration iterative algorithm}  converges for a sufficiently large value of the comparison medium values  $\mathcal{A}_0$, but the rate of convergence depends on the values \cite{moulinec_2018}.  Figure \ref{fig:comp_med} shows the number of iterations required for convergence for various values of the comparison medium for both the model and random geometries.  We see that the best value is the smallest value for which there is convergence.  We use $\mathcal{A}_0 = 0.55, \mathcal{B}_0=1$ in what follows.

\paragraph{Verification.} 
\begin{figure}[t]
\centering
    \includegraphics[width=6in]{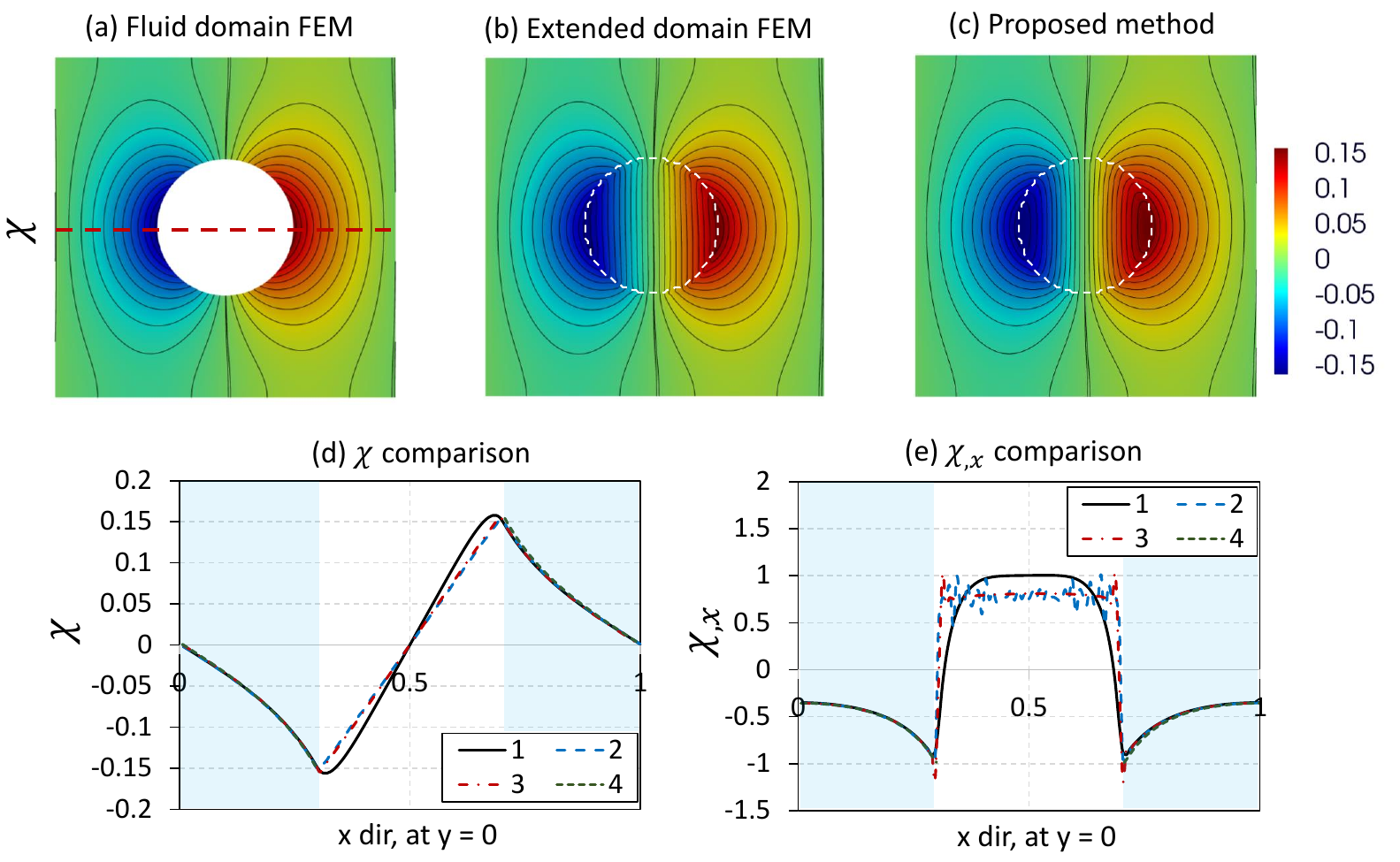}
    \caption{Verification of solute transport in the model geometry.  (a-c)  The concentration profiles computed using the three methods with $\bfg_p= [1 \ 0]$, $\bfg_\chi = [1 \ 0]$, $\eta = 0.01$, $\epsilon=10^{-5}$  (The white dashed line indicates the solid/fluid interface)  (d) Concentration and (e) horizontal gradient of the concentration along the line $y=0.5$ indicated by the dashed line in (a) using the three methods as well as exact FFT. Plotted lines 1, 2, 3, and 4 are associated with approximate FFT, exact FFT, extended domain FEM, and fluid domain FEM, respectively. }
    \label{fig:model_sol}
\end{figure}

We now turn to the verification of the numerical method.  We solve the solute problem with three different approaches.
\begin{itemize}
\item {\it Fluid domain FEM.}  We solve the original transport equation (\ref{eq: unit cell})$_{3}$ on the fluid domain $Y_p$ subject to periodic boundary conditions on the fluid boundary of the unit cell $\partial Y_p \cap \partial Y$ and to Neumann boundary condition (\ref{eq: unit cell})$_{4}$  on the solid-fluid interface $Y_p \cap Y_s$.  We use an unstructured mesh to resolve the equations on $Y_p$.  
\item {\it Extended domain FEM.}  We solve the iterative formulation introduced in Section  \ref{sec:concentration iterative algorithm} on the extended domain by using a finite element discretization of (\ref{eq:iter}) on the unit cell $Y$.
\item {\it Proposed method.} We solve the iterative formulation introduced in Section  \ref{sec:concentration iterative algorithm} on the extended domain using  FFT and GPUs as described in Section \ref{sec:gpu}.
\end{itemize}

Figure \ref{fig:model_sol} (a-c) shows the solute concentration computed using the three methods for the model geometry.  The fluid domain FEM used an unstructured mesh consisting of $100856$ triangular elements, while the other two methods used a structured discretization with the resolution of $128\times 128$.  We observe that the constrained formulation in both the extended domain FEM and proposed FFT-based method results in $(\nabla \chi =  \bff = [1 \ 0]$) inside the solid region as required, and the results of the three methods agree closely.  

Figure \ref{fig:model_sol}(d,e) shows the concentration and horizontal gradient of the concentration along the center line.  The FFT in the proposed method is computed in two ways: exactly, and approximately according to (\ref{eq:fft_app}).  We observe that the exact FFT leads to spurious oscillations that are absent in the two FEM methods.  The approximate FFT filters these away but still provides the correct values in the physically meaningful fluid domain.  The corresponding results for the random geometry are provided in Supplementary Materials, Figure \ref{fig:random_osc1}.  

\begin{figure}
\centering
    \includegraphics[width=4.5in]{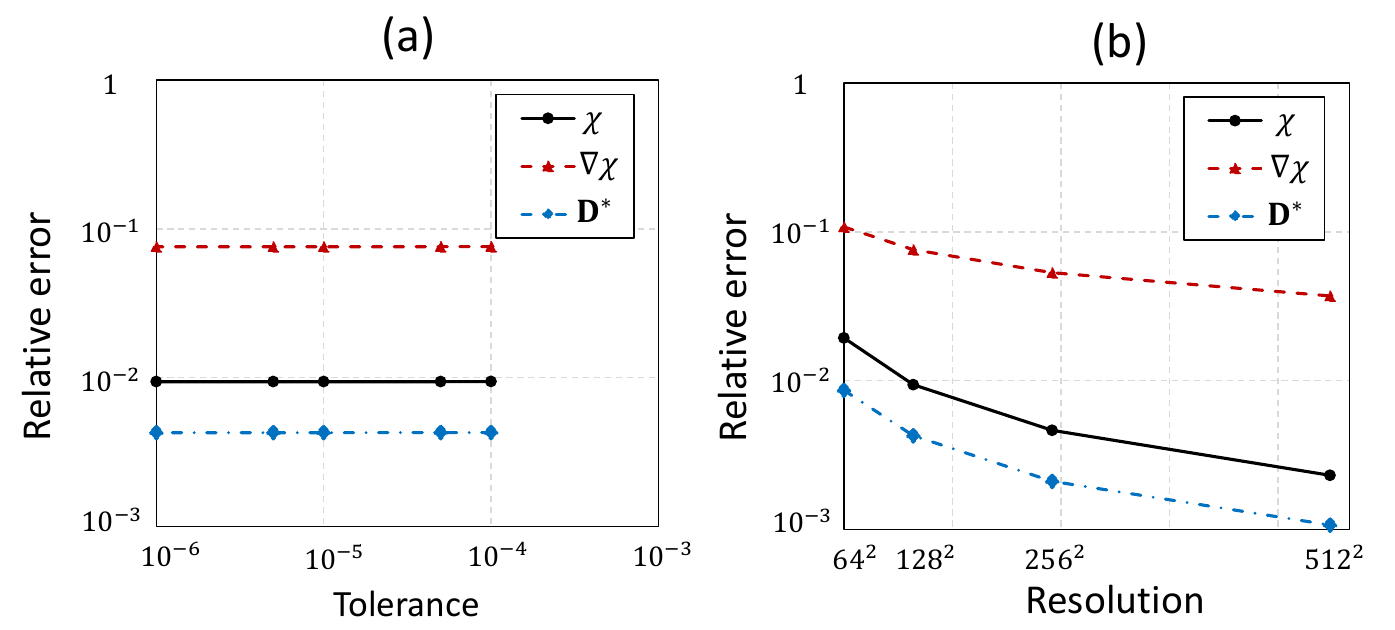}
    \caption{The relative errors in concentration fields and the diffusivity computed using the proposed method relative to the fluid domain FEM in the model geometry for various  (a) tolerances with resolution $128\times 128$ and (b) resolutions with $\epsilon=10^{-5}$. $\bfg_\chi = [1 ~ 0]$.} 
    \label{fig:model_err_chi}
\end{figure}

The relative errors in concentration fields and the diffusivity for various tolerances and resolution for the model geometry are shown in Figure \ref{fig:model_err_chi}.

\begin{figure}
 \centering
     \includegraphics[width=5.5in]{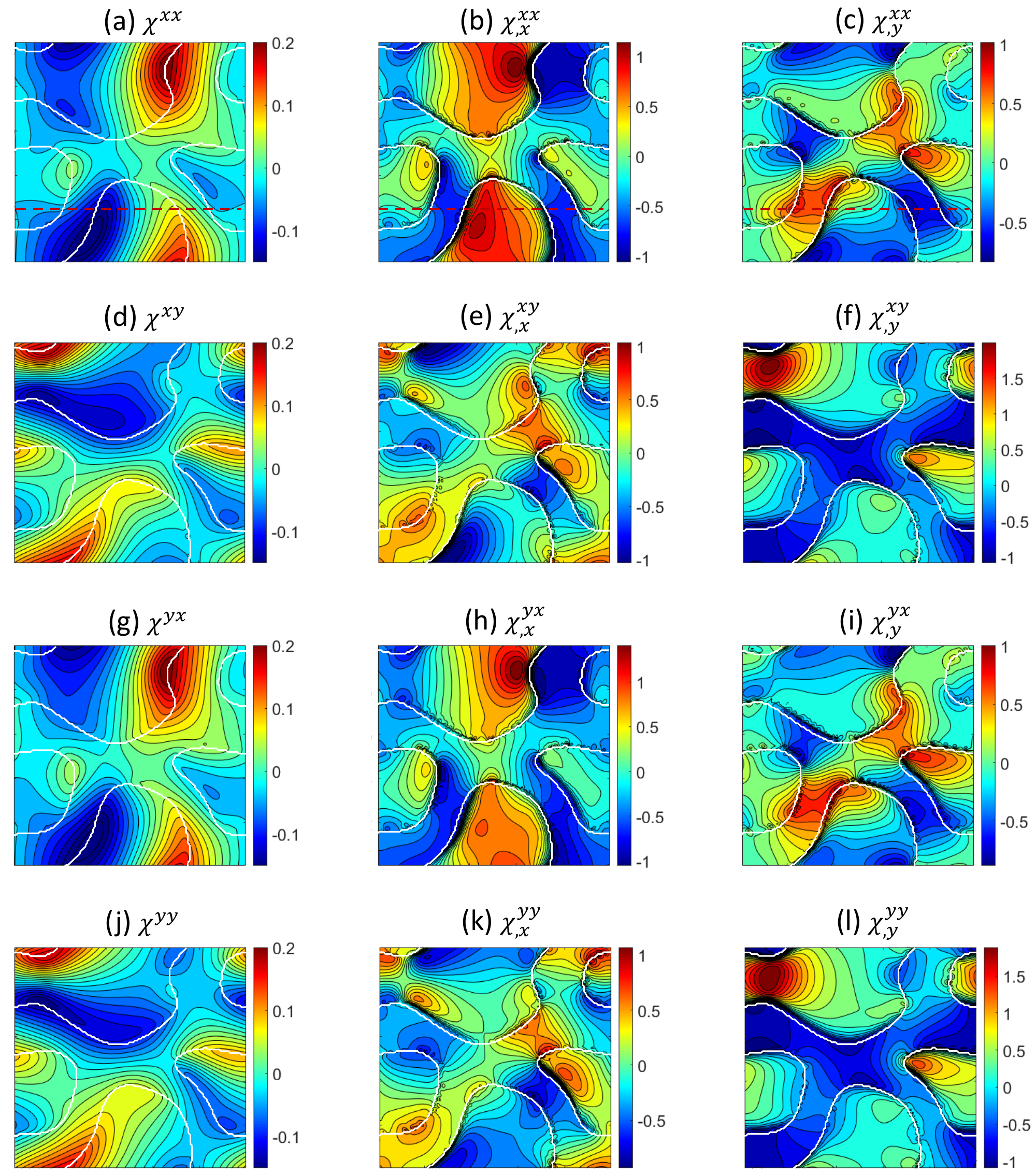}
     \caption{Concentration fields and their gradients for various imposed pressure and composition gradients.  $\chi^{xx}$ corresponds to $\bfg_p = [1 \ 0], \bfg_\chi = [1 \ 0]$, $\chi^{xy}$ corresponds to $\bfg_p = [1 \ 0], \bfg_\chi = [0 \ 1]$, $\chi^{yx}$ corresponds to $\bfg_p = [0 \ 1], \bfg_\chi = [1 \ 0]$, and $\chi^{yy}$ corresponds to $\bfg_p = [0 \ 1], \bfg_\chi = [0 \ 1]$.       $\tt{Pe = 50}$.}
     \label{fig:random_conc}
 \end{figure}
%
%

We conclude by noting the anisotropy and asymmetry in the concentration due to the geometry, Figure \ref{fig:random_conc}.   This figure shows the concentration and its gradients for various imposed pressure and composition gradients.  The first and the third rows (also the second and the fourth rows) have the same imposed composition gradient, but different distribution because of the difference in flow induced by the different pressure gradients.   

\subsection{Computational perforrmance}

\paragraph{Computational cost}
We compare the computational cost of the three methods introduced above for solving the fluid flow as well as the solute transport. The two FEM-based methods are solved on a single core of an Intel Skylake CPU with a clock speed of 2.1 GHz, while the FFT-based method is solved both on the same CPU and on a NVIDIA P100 GPU with 3584 cores and 1.3 MHz clock speed. 

We find that the wall clock time required for the proposed FFT-based method on a GPU is an order of magnitude faster than all the other methods.  The proposed method implemented on a CPU is comparable to that of the FEM approach on the fluid domain.  The iterative method using full domain resolution is significantly more expensive.

\begin{figure}
\centering
    \includegraphics[width=0.75\textwidth]{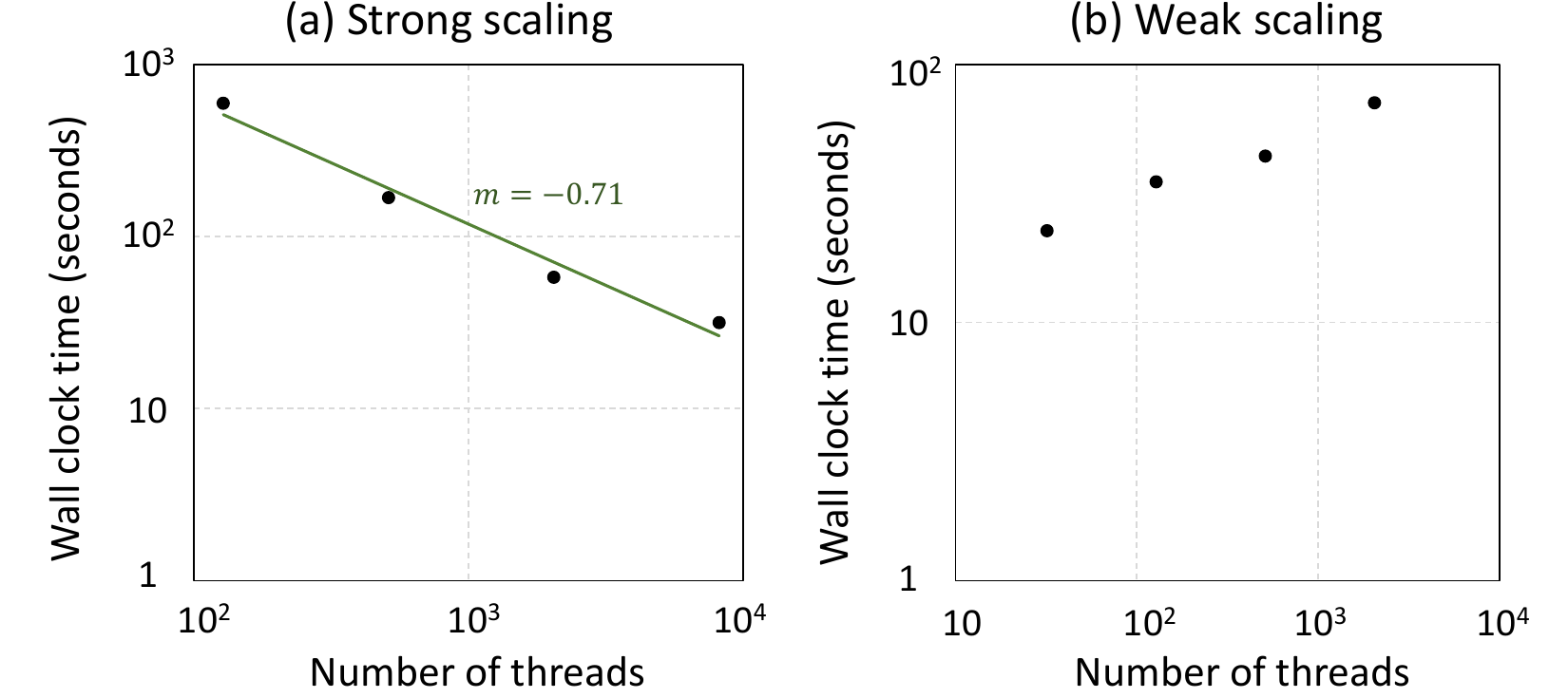}
    \caption{ Performance of the proposed method. (a) Strong scaling and (b) weak scaling.}
    \label{fig:scaling}
\end{figure}

\paragraph{Parallel performance}
Figure \ref{fig:scaling} shows the strong and weak scaling of the proposed FFT-based algorithm as implemented on a GPU.  For the strong scaling shown in Figure \ref{fig:scaling}(a), we consider the problem at a of $1024\times 1024$, and implement it with 128, 512, 2048, and 8192 threads.  We see that the slope on the log-log plot of wall clock time vs. threads is -0.71, demonstrating good parallel performance.  Figure \ref{fig:scaling}(b) shows the weak scaling where we increase the number of threads proportionally to resolution.  We show results for a resolution of $128\times128$, $256\times256$, $512\times512$, and $1024\times1024$ are considered, with 32, 128, 512, and 2048 threads respectively.  We see that the wall clock is almost constant, again demonstrating good parallel performance.  The corresponding results for the random geometry are shown in Supplementary Materials, Figure \ref{fig:scaling-random}.

\begin{table}
    \centering
    \caption{Comparison of computation cost (wall-clock time in seconds)}\label{tab:1}
    \begin{tabular}{ c c } 
    \hline
    Method \quad \quad \quad \quad \quad \quad & Cost \\
    \hline
    FEM (Fluid domain) \quad \quad \quad \quad \quad \quad & $46$ (CPU) \\ 
    FEM (Extended domain*) \quad \quad \quad \quad \quad \quad &  $558$ (CPU) \\
    FEM (Extended domain) \quad \quad \quad \quad \quad \quad &  $10360$ (CPU) \\
    FFT \quad \quad \quad \quad \quad \quad &  $52$ (CPU) \\
    FFT (Parallel) \quad \quad \quad \quad \quad \quad & $2.7$ (GPU) \\
    \hline
    *without considering the auxiliary velocity \\
    \hline
    \end{tabular}
\end{table}

\section{Conclusion}
We develop and demonstrate an efficient computational method to study fluid flow and solute transport problems in periodic porous media.  The work is motivated by the challenges posed by reactive transport in permeable porous media that arise in a number of applications.   We consider periodic porous media where the fluid flow is governed by the Stokes equation, and solute transport is governed by the advection-diffusion equation.   The computational approach proposed here combines two ideas. The first is to use an extended domain and a computational grid that does not conform to the geometry of the porous medium. The second is to formulate the governing equations in a manner that can exploit the massively parallel architecture of GPUs. 

We follow complementary approaches to fluid flow and solute transport.  The Stokes equation governing fluid flow may be written as a variational principle.  So we include incompressibility and the zero velocity in the solid as constraints that are imposed using an augmented Lagrangian.  We then follow the alternating direction method of multipliers (ADMM)  to derive an iterative method.  The advection-diffusion equation governing solute transport can not be written as a variational principle.  Therefore, we introduce a (small) fictitious diffusivity in the solid region and solve the problem using a Lippmann-Schwinger approximation.  In both cases, the iterative method consists of local, and hence trivially parallelizable, equations and a Poisson's equation, which also can be effectively parallelized using fast Fourier transforms.  We exploit this structure to implement the methods on accelerators (GPUs).  We discuss the convergence of these methods and verify them against solutions obtained using the finite element approach.  We then demonstrate the efficacy of the approach by comparing the performance against other methods and also demonstrating strong and weak parallel scaling. 

This work is motivated by multiscale modeling where the unit cell problem studied here is used to infer macroscopic properties of the porous medium over scales that are orders of magnitude larger than the pore size.   This is especially useful in situations where the flow results in the deposition on or dissolution of the porous solid.  In such situations, the porous microgeometry, and consequently the flow and transport at the microscale as well as the effective transport property at the macroscale, can change over time.  Therefore, one has to repeatedly update the fluid flow and solute transport repeatedly as the porous microgeometry changes.  The proposed approach addresses this challenge effectively.   First, the numerical discretization is independent of the geometry in the proposed approach, and one can repeat the calculation with the same discretization as the geometry changes.  Second, the speed and efficiency of the method allow us to do numerous calculations at a modest cost.

We use the methods developed here in Karimi and Bhattacharya \cite{karimi2023learning} in a multiscale setting where one has deposition and dissolution.  
Specifically, we solve the unit cell problem with various histories of concentration and pressure and use it as the data to train a recurrent neural operator that describes the history-dependent change of macroscopic properties.  

In this work, we have implemented the framework to two dimensions to keep computational cost modest.  However, the approach generalizes naturally to three dimensions, and this is the topic of current research.  Finally, we address Stokes flow and advection reaction equation in this work.  However, the approach can be generalized to a wide variety of problems in porous media, including poroelasticity and phase transformations.

\section*{Acknowledgments}
We gratefully acknowledge the financial support of the Resnick Sustainability Institute at the California Institute of Technology.  KB also acknowledges the support of the Army Research Office through grant number W911NF-22-1-0269.  The simulations reported here were conducted on the Resnick High Performance Computing Cluster at the California Institute of Technology.


\newpage
\begin{center}
Karimi and Bhattacharya\\
A fast-Fourier transform method for reactive flow in porous media\\
\vspace{0.1in}
{\Large \bf Supplementary Materials}\\
\end{center}
\renewcommand{\thepage}{S\arabic{page}}
\renewcommand{\thefigure}{S\arabic{figure}}
\setcounter{page}{1}
\setcounter{figure}{0}

\begin{figure}[h]
\centering
    \includegraphics[width=0.8\textwidth]{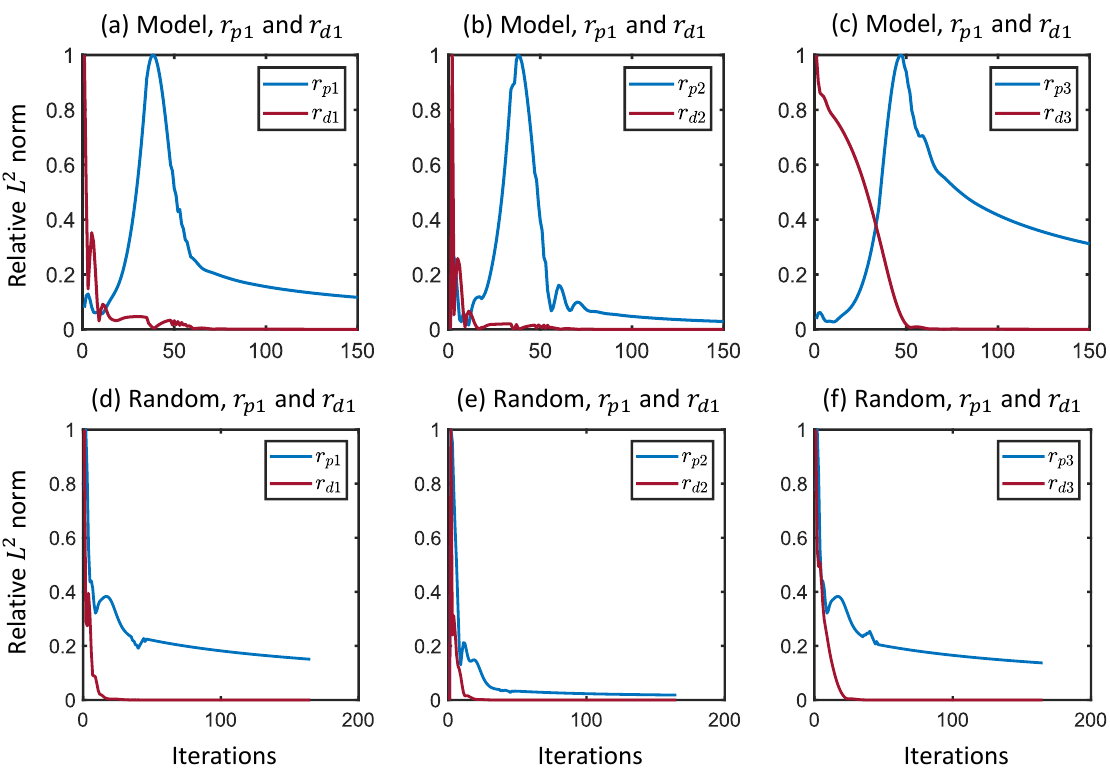}
    \caption{Primal and dual errors (relative $L_2$ norm normalized by the highest value during the iteration) of three constraints associated with the three constraints the model geometry (a-c) and random geometry (a-c).  $\bfg_p = [1,0]$ in these simulations.}
    \label{fig:prim_dual}
\end{figure}

\begin{figure}
\centering
    \includegraphics[width=0.85\textwidth]{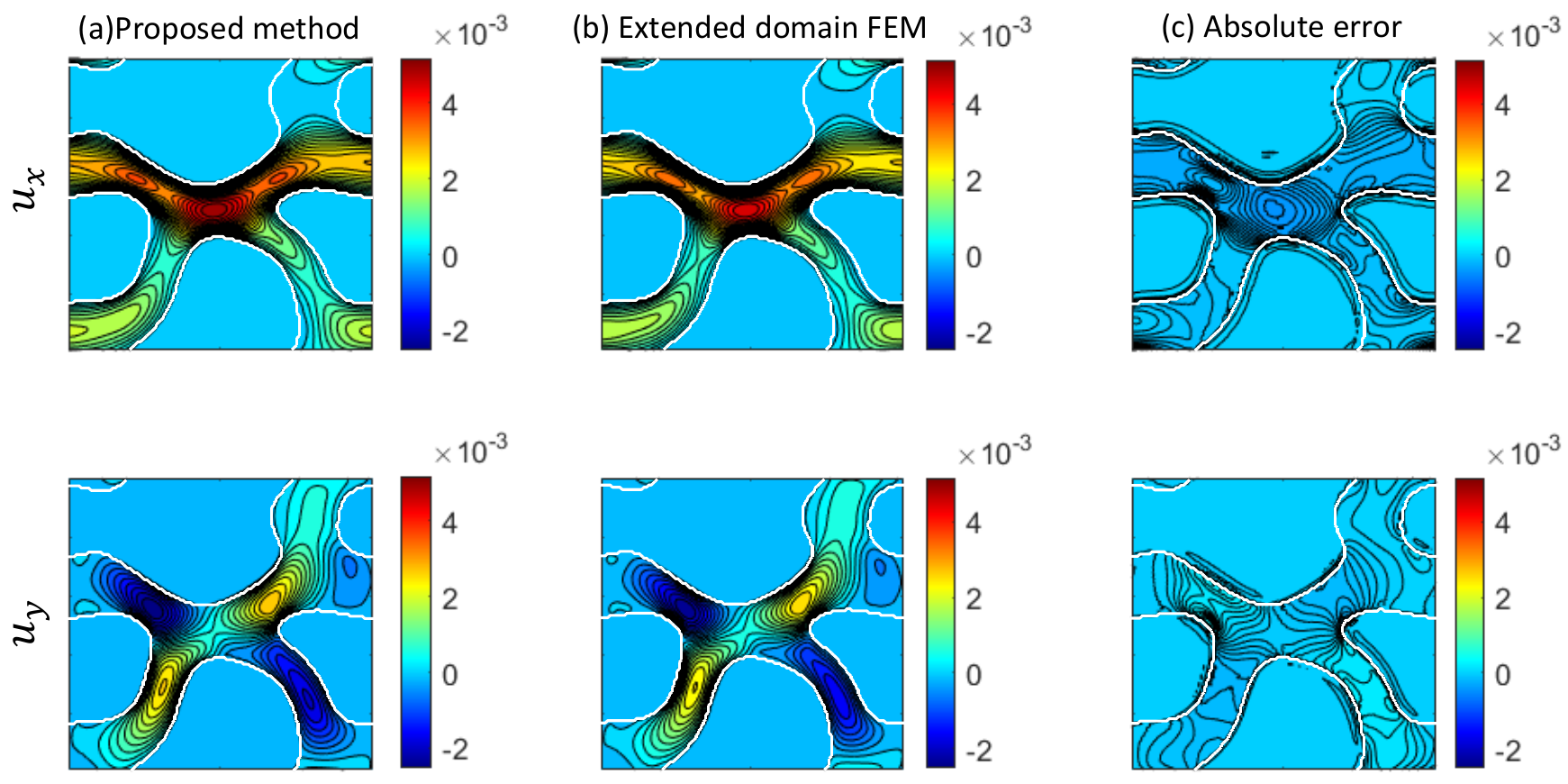}
    \caption{Comparison of velocity profiles computed using the extended domain FEM and proposed method in the random geoemetry with $\bfg_p = [1 ~ 0]$ and $\epsilon = 10^{-5}$ (the white dashed line indicates the solid/fluid interface).}
    \label{fig:random_vel1}
\end{figure}

\begin{figure}
\centering
    \includegraphics[width=0.85\textwidth]{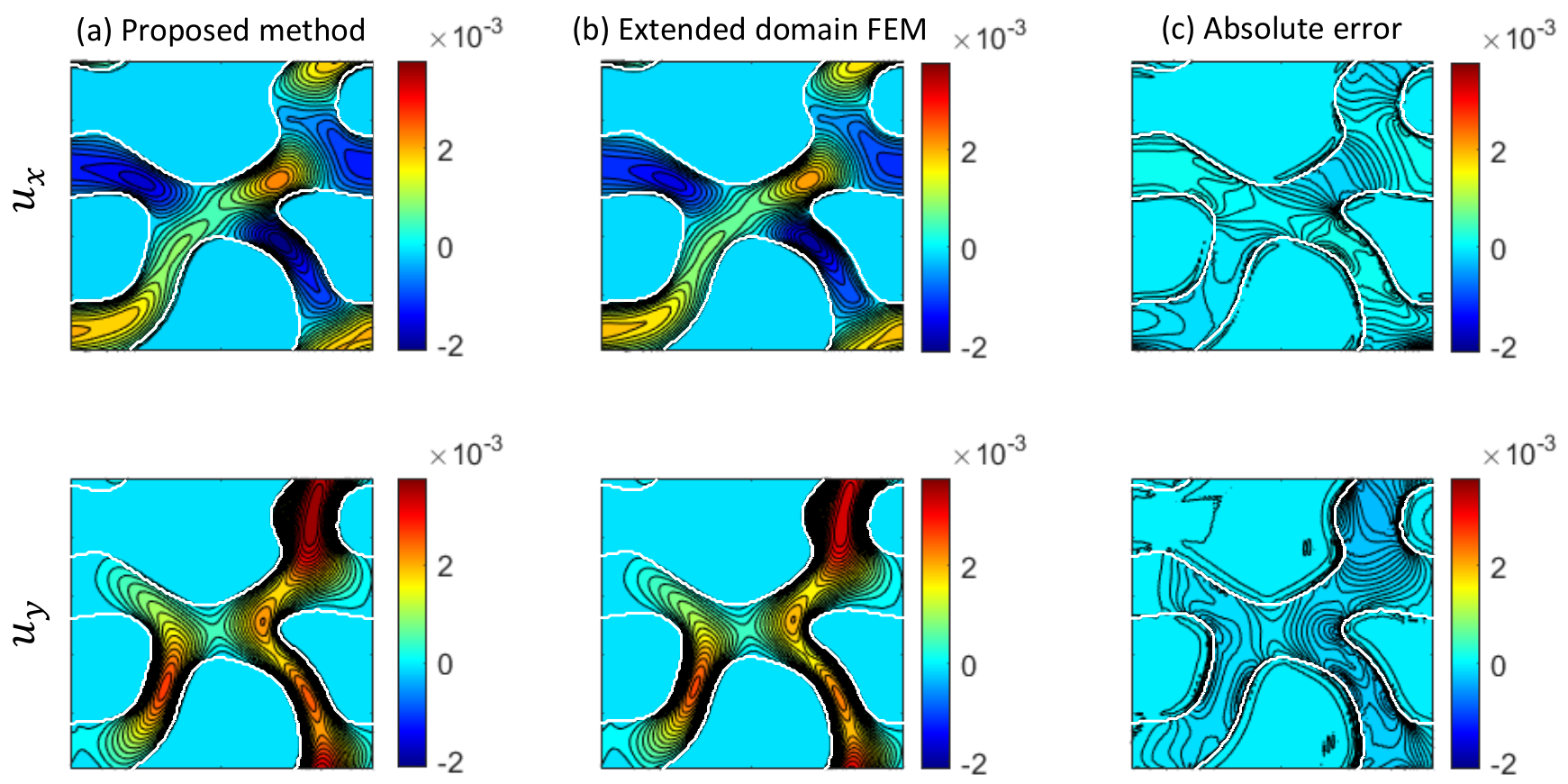}
    \caption{Comparison of velocity profiles computed using the extended domain FEM and proposed method  in the random geoemetry with $\bfg_p = [0 ~ 1]$ and $\epsilon = 10^{-5}$ (the white dashed line indicates the solid/fluid interface).}
    \label{fig:random_vel2}
\end{figure}

\begin{figure}
\centering
    \includegraphics[width=0.8\textwidth]{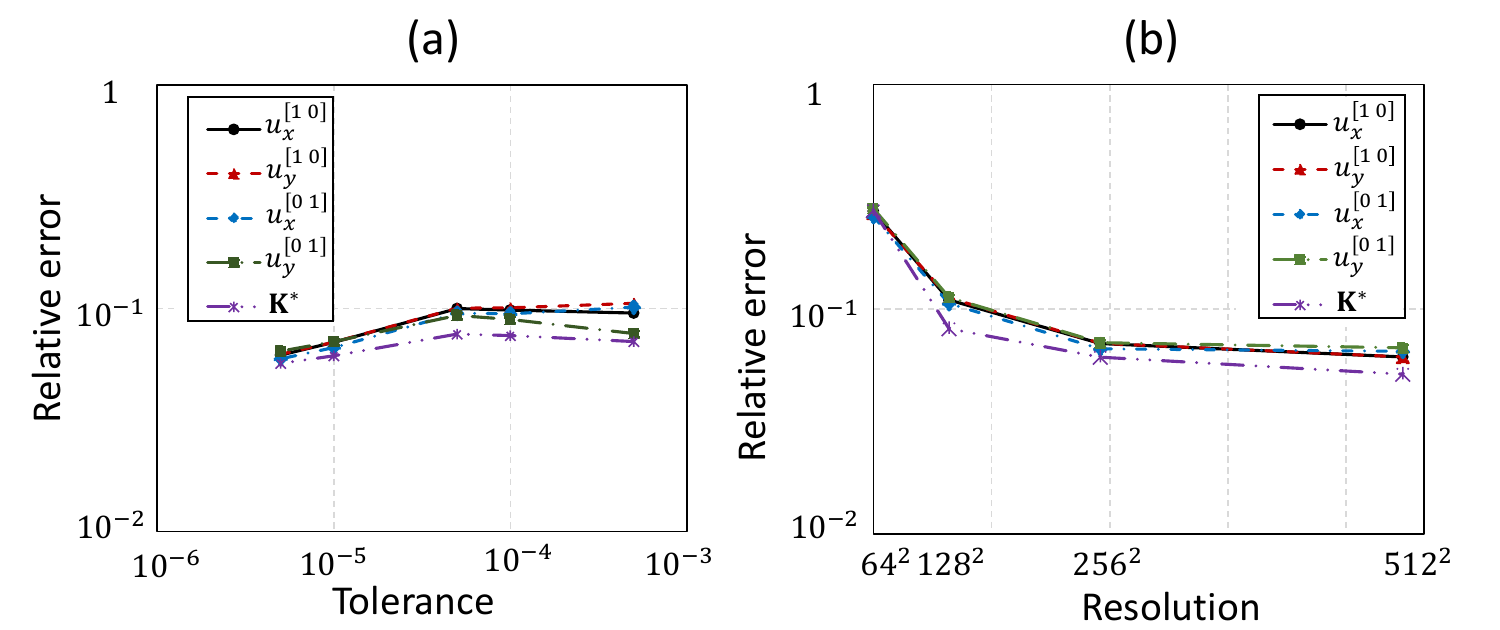}
    \caption{The relative errors in average quantities and velocity in the random geometry computed using the proposed method compared to the extended domain FEM for various (a) tolerances with resolution $256\times 256$ and (b) resolutions with tolerance $\epsilon = 10^{-5}$. }
    \label{fig:random_err}
\end{figure}


\begin{figure}
\centering
    \includegraphics[width=6in]{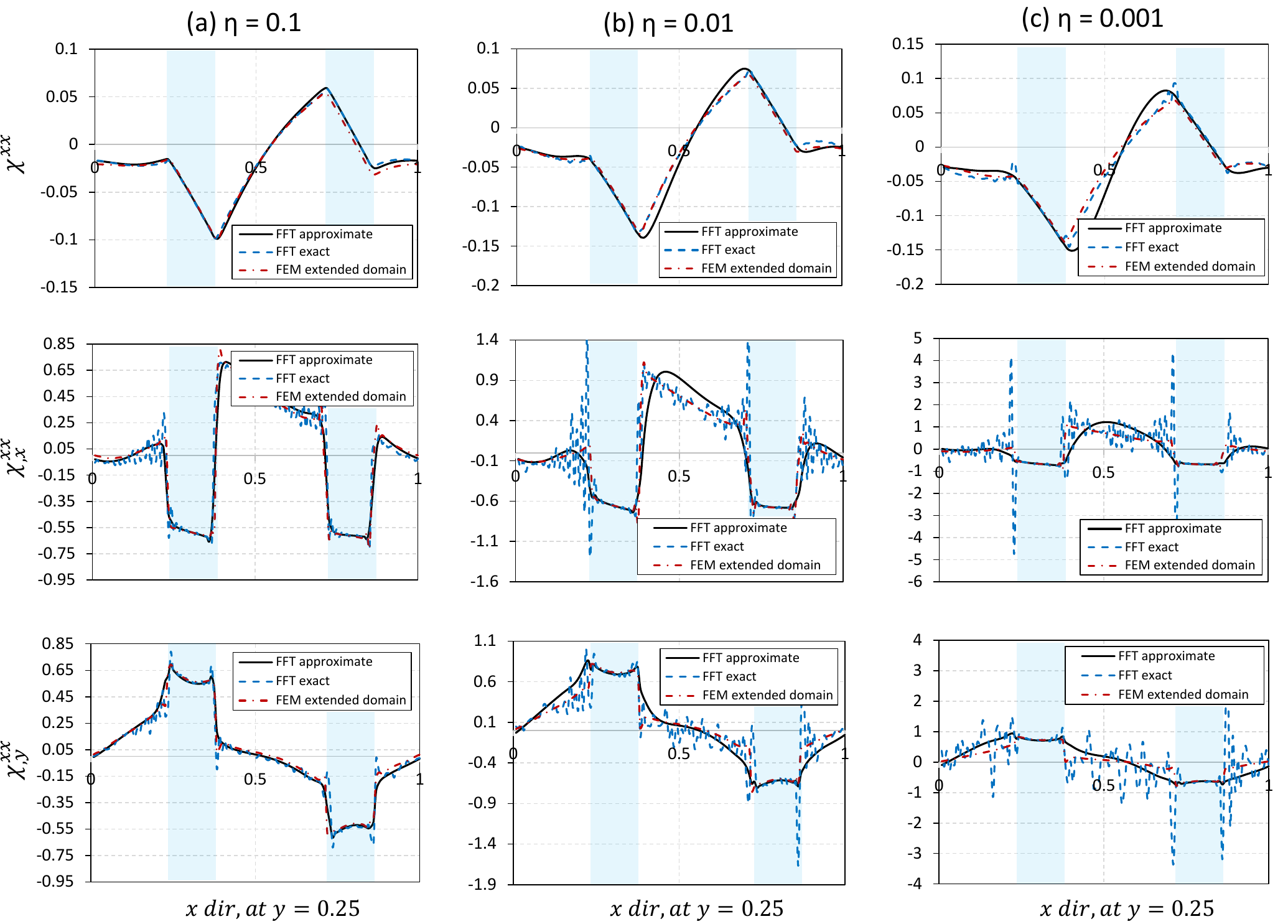}
    \caption{Comparison of scalar concentration field $\chi^{xx}$ and concentration gradients in the $x$, and $y$ directions ($\chi^{xx}_{,x}, ~ \chi^{xx}_{,y}$), using exact FFT, approximate FFT and FEM approaches, considering different values of $\eta$. The results are extracted in the $x$ direction, along the dashed line shown in figure \ref{fig:random_conc}, considering $\bfg_p = [1 \ 0]$, and $\bfg_\chi = [1 \ 0]$. $\tt{Pe = 50}$. }
    \label{fig:random_osc1}
\end{figure}

\begin{figure}
\centering
    \includegraphics[width=6in]{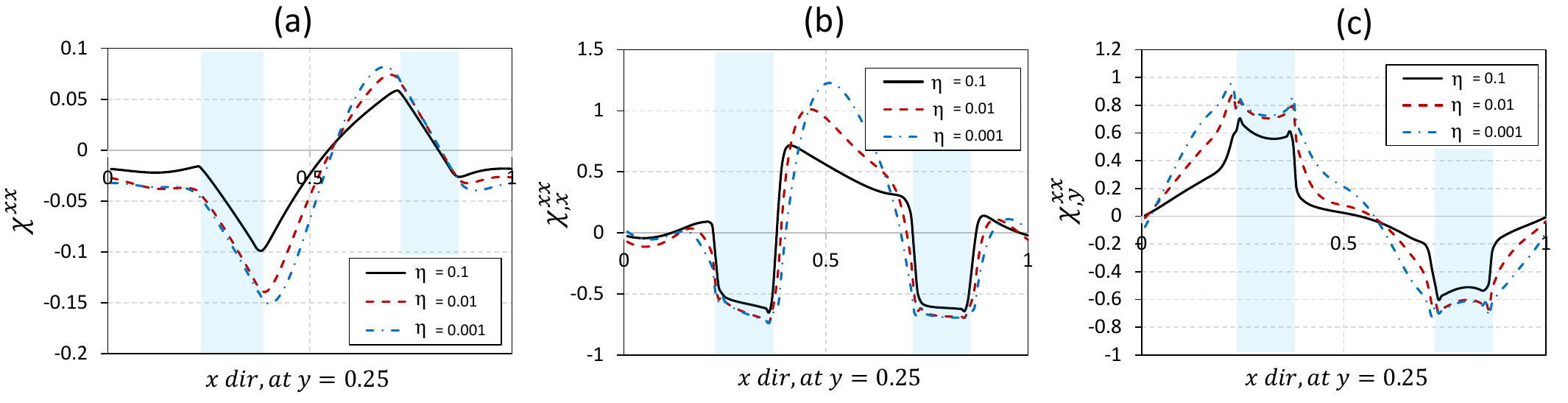}
    \caption{Variation of scalar concentration field $\chi^{xx}$ and concentration gradients in the $x$, and $y$ directions ($\chi^{xx}_{,x}, ~ \chi^{xx}_{,y}$), considering different values of $\eta$, using the approximate FFT method. The results are extracted in the $x$ direction, along the dashed line shown in figure \ref{fig:random_conc}, considering $\bfg_p = [1 \ 0]$, and $\bfg_\chi = [1 \ 0]$. $\tt{Pe = 50}$.}
    \label{fig:random_osc2}
\end{figure}

\begin{figure}
\centering
    \includegraphics[width=0.75\textwidth]{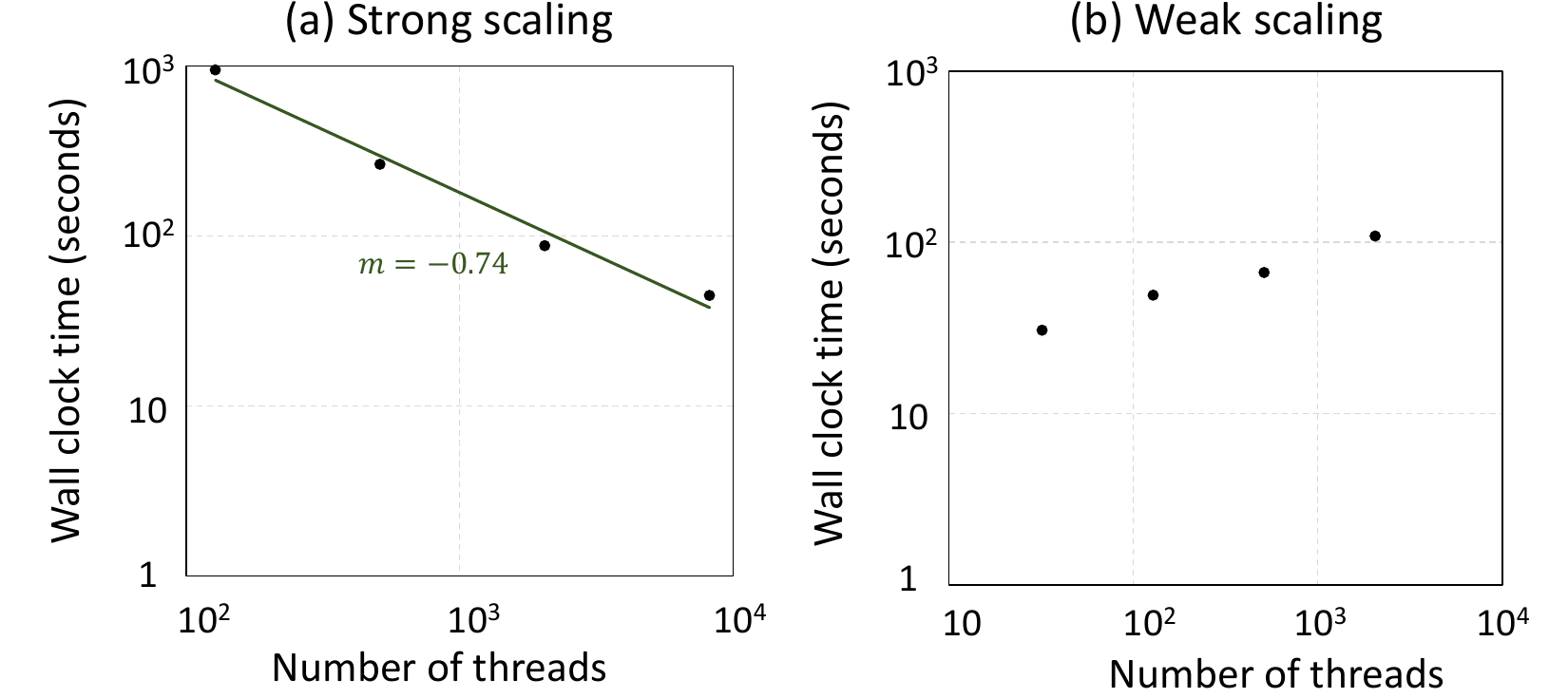}
    \caption{ Performance of the proposed FFT-based method calculated for random geometry. (a) Strong scaling considering a resolution of $1024\times 1024$, and (b) weak scaling considering resolutions $128\times 128$, $256\times 256$, $512\times 512$, and $1024\times 1024$, correspond to thread numbers $32$, $128$, $512$, and $2048$, respectively. $\epsilon =10^{-6}$. }
    \label{fig:scaling-random}
\end{figure}

\end{document}